\documentclass[a4paper]{amsart}

\usepackage{amsmath}
\usepackage{amsfonts}
\usepackage{amssymb}
\usepackage{latexsym}
\usepackage{epsfig}

\def\co{\colon\thinspace}
\newcommand{\begriff}[1]{\textbf{#1}}

\renewcommand{\d}{\partial}

\renewcommand{\t}{\mathcal T}
\newcommand{\CR}{\mathop{Cr}}
\renewcommand{\c}{\mathcal C}

\newtheorem{theorem}{Theorem}

\theoremstyle{definition}
\newtheorem{definition}{Definition}

\theoremstyle{remark}

\begin{document}
\title{Crossing number  of links formed by edges of a triangulation}  
\author{Simon A. King}
\address{Simon A. King\\
Department of Mathematics\\
Darmstadt University of Technology\\
Schlossgartenstr.~7\\
64289 Darmstadt\\
  Germany}
\email{king@mathematik.tu-darmstadt.de}
\begin{abstract}
  We study the crossing number of links that are formed by edges of a
  triangulation $\t$ of $S^3$ with $n$ tetrahedra. 
  We show that the crossing number is bounded from above by an
  exponential function of $n^2$. In general, this bound can not be
  replaced by a subexponential bound. 
  However, if $\t$ is polytopal (resp. shellable) then there is a
  quadratic (resp. biquadratic) upper bound in $n$ for the crossing
  number. 
  In our proof, we use a numerical invariant $p(\t)$, called
  \emph{polytopality}, that we have introduced in~\cite{king2}. 
  \\
  \textbf{Primary:}        57M25, 
                                              57Q15  
  \\
  \textbf{Secondary:}  52C45,  
                                               52B22 
\end{abstract}

\maketitle

\section{Introduction}
\label{sec:intro}

The aim of this paper is to study properties of knots and links that are
formed by edges of a triangulation of the $3$--dimensional sphere
$S^3$. 
For any $n\in\mathbb N$, let $\mathbb T_n$ denote the set of
triangulations of $S^3$ with $n$ tetrahedra. The number of
simplicial complexes with $n$ tetrahedra is finite, and the number of
simple edge paths in a simplicial complex is finite as well. 
Therefore, there are only finitely many equivalence classes of links in
$S^3$ occuring in the $1$--skeleton $\t^1$ of some $\t\in\mathbb T_n$.
Thus, if $f(\cdot)$ is any numerical link invariant (e.g., bridge
number, crossing number, $k$--th coefficient of the Jones polynomial)
then 
$$F(n) = \sup_{\t\in\mathbb T_n,\;L\subset \t^1} f(L)$$  
is a finite number for all $n\in \mathbb N$.

In~\cite{king1} and~\cite{king2}, we obtained exponential upper and
lower bounds for $F(n)$ when $f(L)$ is the bridge number of $L$. 
However, there are infinitely many non-equivalent links whose bridge
number is smaller than $F(n)$, for instance, two-bridge knots. 
The aim of this paper is to  prove an upper bound for $F(n)$ when $f(L)$ is
the crossing number of $L$. Such an estimate is substantially stronger than an
estimate for the bridge number, in the sense that there are only
finitely many equivalence classes of links with a crossing number smaller
than $F(n)$. 

We recall the notions of crossing number  and bridge number of a link
$L\subset S^3$.  
The \begriff{crossing number} $\CR(L)$ of $L$ is the minimal number of
crossings in a link diagram representing $L$. This is a natural measure of 
complexity of links, and was already used in the knot tables of
Tait and Little in 1900.
The bridge number has been introduced by H.~Schubert~\cite{schubert} in
1954. It can be defined as follows.
Let $I=[0,1]$ denote the unit interval. 
Let  $H\co S^2\times I\to S^3$ be a smooth embedding with $L\subset
H(S^2\times I)$. Let $\pi_I\co S^2\times I\to I$ denote the projection
to the second factor.
On gets a map $h\co L\to S^2$ by $h(p) = \pi_I\circ
H^{-1}(p)$, for $p\in L$.  
We call $p\in L$ a \begriff{critical point} of $H$ with respect to $L$
if $H(S^2\times \{h(p)\})$ is not transversal to $L$ in $p$. Thus,
generically a critical point is an isolated local maximum or minimum of
$h$.  We denote the number of critical points by $c(H,L)$.
Now, the \begriff{bridge number} of $L$ is
$$ b(L) = \frac 12 \min_H c(H,L),$$ 
where the minimum is taken over all embeddings $H\co S^2\times I\to S^3$
with $L\subset H(S^2\times I)$.
The factor $\frac 12$ is needed to make this definition consistent with
Schubert's original definition~\cite{schubert}.

Link invariants also arise in the study of polytopal and
shellable triangulations of $S^3$, as explained in the next
paragraphs. Recall that a triangulation of the 
$d$--dimensional sphere $S^d$  is \begriff{polytopal}, if it is
isomorphic to the boundary complex of a convex $(d+1)$--polytope. A
triangulation of $S^d$ is \begriff{shellable}, if there is an order
$\sigma_1,\dots, \sigma_{n}$ of its $d$--simplices so that
$\bigcup_{i=1}^k \overline{\sigma_i}$ is a $d$-dimensional closed
ball, for $k=1,\dots,n-1$. For a more general definition of shellable cellular complexes,
see~\cite{ziegler}. 
Any polytopal triangulation is shellable~\cite{bruggesser}. There are
shellable triangulations of $S^3$ that are not polytopal, and there are
triangulations of $S^3$ that are not shellable~\cite{ziegler}.

Let $\t\in \mathbb T_n$ be a triangulation of $S^3$ with $n$ tetrahedra, and let
$K\subset \t^1$ be a {knot} formed by $k$ edges of $\t$.  
It follows from the work of Lickorish~\cite{lickorish} that 
$b(K)\le k$ under the assumption that $\t$ is shellable. 
Armentrout~\cite{armentrout} obtained $b(K)\le \frac 12 k$
under the assumption that the \emph{dual} of $\t$ 
is shellable. 
Recently Ehrenborg and Hachimori~\cite{ehrenborghachimori} obtained
the sharp bounds $b(K)\le \frac 12 k$ if $ T$ is shellable and $b(K)\le \frac 13 k$ if
$\t$ is vertex decomposable (we will not define this notion here).
In the same way, one can show that if $L\subset  \t^1$ is a link (but not
necessarily a knot) then $b(L)$ is bounded from above by a linear
function of $n$, provided $\t$ satisfies one of the mentioned
assumptions.  
These assumptions are in fact very strong,  since in general
there is no subexponential upper bound for $b(L)$ in terms of $n$,
see~\cite{king2}. Without any geometric assumption on $\t$, we found 
$b(L) < 2^{190n^2}$ by a complexity  analysis of the
Rubinstein--Thompson algorithm for the recognition of $S^3$, see~\cite{king1}.    

By the following main theorem of this paper, the crossing number of a
link in $\t^1$ is also sensitive for geometric properties of $\t$.
\begin{theorem}\label{thm:crossing}
  Let $\t$ be a triangulation of $S^3$ with $n$ tetrahedra, and let
  $L\subset  \t^1$ be a link.  
  \begin{enumerate}
  \item If $\t$ is polytopal then $\CR(L)< 4n^2$.
  \item If $\t$ or its dual cellular decomposition is shellable then
    $\CR(L)< 10^9 n^4$. 
  \item In general, $\CR(L) < 2^{810 n^2}$.
  \end{enumerate}
\end{theorem}
So far as known to the author, these are the first upper bounds for 
$\CR(L)$ in terms of $\t$. Certainly these bounds are not optimal, we
did not try to prove sharp bounds.
However, since $\CR(L)\ge b(L)$, our results in~\cite{king2} imply that the
general bound in the third part of Theorem~\ref{thm:crossing}  can not
be replaced by a subexponential bound. Thus, the assumptions in the first two
parts of Theorem~\ref{thm:crossing} can not be removed.

We outline the proof of Theorem~\ref{thm:crossing}.
Let $\t\in\mathbb T_n$ and let $L\subset  \t^1$.
If $\t$ is polytopal, then one can isotope its $2$--skeleton in $\mathbb
R^3\subset S^3$ so that all edges become straight line segments. 
Thus, an orthogonal projection yields a link diagram for $L$ with at
most one crossing for each pair of edges. This yields the first part of
Theorem~\ref{thm:crossing}. 

The idea for the proof of the second and third part of
Theorem~\ref{thm:crossing} is to transform the triangulation $\t$ into a 
polytopal triangulation $\tilde\t$ by a finite number of local changes
(e.g., stellar subdivisions), so that $\tilde\t^1$ contains a copy
$\tilde L$ of $L$, where the number of edges forming $\tilde L$  is
controlled in terms of the number of local changes. This reduces the
problem to the first part of Theorem~\ref{thm:crossing}.  
For this purpose, we use the  notion of \emph{polytopality} that we have
introduced in~\cite{king2}.   
The polytopality $p(\t)$ is a numerical invariant of a triangulation 
$\t$ of $S^3$. Its definition is inspired by the bridge number of links,
see Section~\ref{sec:polytop}.  
In~\cite{king2}, we have shown how $\tilde\t$ can be constructed from 
$\t$, so that the number of local changes is bounded from above
by a quadratic function of $p(\t)$. 
Since $\tilde\t$ is polytopal, we obtain a biquadratic upper  bound for 
$\CR(\tilde L) = \CR(L)$ in terms of $p(\t)$.  
We have shown~\cite{king2} that $p(\t)\le 7n$ if $\t$ or its dual is
shellable, and $p(\t)< 2^{200 n^2}$ in general. This yields the
estimates claimed in the second and third part of 
Theorem~\ref{thm:crossing}.

The rest of this paper is organized as follows. In
Section~\ref{sec:polytop} we define the polytopality $p(\t)$ of a
triangulation $\t$ of $S^3$ and recall from~\cite{king2} how to
construct from $\t$ a polytopal triangulation by a series
of local changes whose length is bounded in terms of $p(\t)$.
In Section~\ref{sec:proof}, we start with a proof of the first part of
Theorem~\ref{thm:crossing}. Then, we formulate  and prove a bound for
the crossing number of a link $L\subset \t^1$ in terms of $p(\t)$. This
bound together with estimates for $p(\t)$ from~\cite{king2} immediately
yields the second and third part of Theorem~\ref{thm:crossing}.

\section{Polytopality of triangulations}
\label{sec:polytop}

Let $\t$ be a triangulation of $S^3$. The aim of this section is to expose a numerical
invariant $p(\t)$, called polytopality. Although its definition is purely
topological, $p(\t)$ turns out to be a measure for the {geometric}
complexity of $\t$. 
 
We recall the definition of $p(\t)$ from~\cite{king2}.
Let $\mathcal C$ be the dual cellular decomposition of $\t$, with
$1$--skeleton $\mathcal C^1$.
Let $H\co S^2\times I\to S^3$ be a smooth embedding with $\c^1\subset
H(S^2\times I)$. Let $\pi_I\co S^2\times I\to I$ denote the projection
to the second factor. 
This gives rise to a map $h\co \c^1\to S^2$ by $h(p) = \pi_I\circ
H^{-1}(p)$, for $p\in \c^1$.  
We call $p\in \c^1$ a \begriff{critical point} of $H$ with respect to $\c^1$
if $H(S^2\times \{h(p)\})$ is not transversal to $\c^1$ in $p$. Thus,
generically a critical point is a vertex of $\c^1$ or an isolated local
maximum or minimum of $h$ in the interior of an edge of $\c$.  
We denote the number of critical points of $H$ with respect to $\c^1$ by
$c(H,\c^1)$ and define 
$$ p(\t) = \min_H c(H,\c^1),$$
the \begriff{polytopality} of $\t$, where the minimum is taken over all
embeddings $H\co S^2\times I\to S^3$ with $\c^1\subset H(S^2\times I)$. 

The definition of polytopality is very similar to the definition of bridge
number in Section~\ref{sec:intro}. One can think of $p(\t)$ as a bridge
number of the spatial graph $\c^1\subset S^3$. 
We remark that already the \emph{abstract} graph $\c^1$ contains much
information about $\t$. For instance, if $\t$ is polytopal than it is
already 
determined by the abstract graph $\c^1$, see~\cite{ziegler}.
In~\cite{king2}, we obtained the following estimates.
\begin{theorem}[see~\cite{king2}]\label{thm:polybounds}
  Let $\t\in \mathbb T_n$.  
  \begin{enumerate}
  \item If $\t$ is polytopal then $p(\t)=n$.
  \item If $\t$ or its dual cellular decomposition is shellable then
    $p(\t)\le 7n$. 
  \item In general, $n\le p(\t) < 2^{200 n^2}$. The upper bound can not be
    replaced by a subexponential bound.\qed
  \end{enumerate}
\end{theorem}

It is not known to the author whether there is a non-polytopal
triangulation $\t\in \mathbb T_n$ with $p(\t)=n$.
Theorem~\ref{thm:polybounds} shows that $p(\t)$ is sensitive
for geometric properties of $\t$. However, the connection to geometry is
even stronger, since $p(\t)$ measures to what extend $\t$ fails to be
polytopal, as explained in the rest of this section. 

\begin{definition}
  Let $\t_1$ and $\t_2$ be triangulations of $S^3$, and let  $e$ be an
  edge of $\t_2$ with $\d e=\{a,b\}$. Suppose that $\t_1$ is obtained from
  $\t_2$ by removing the open star of $e$ and identifying the join
  $a*\sigma$  with $b*\sigma$ for any simplex $\sigma$ in the link of
  $e$. 
  Then $\t_2$ is the result of an \begriff{expansion} of $\t_1$ along
  $e$.  
\end{definition}
Since an expansion increases the number of vertices by
one and the number of simplicial complexes with a given number of
vertices is finite, it is easy to see that the number of possible
expansions of $\t_1$ is finite up to isotopy. 
Any stellar subdivision of a simplex of $\t$ is an expansion. 

It is a consequence of the ``Hauptvermutung'' for $3$--dimensional
manifolds~\cite{moise} that any triangulation $\t$ of $S^3$ can be
turned into a polytopal triangulation by a series of expansions. We
define $d(\t)$ as the length of a shortest series of expansions that
transforms $\t$ into a polytopal triangulation. Our next theorem
shows that $d(\t)$, defined in terms of discrete geometry, is closely
related to $p(\t)$, defined in terms of topology.
\begin{theorem}[see~\cite{king2}]\label{thm:makepoly}
  For any $\t\in\mathbb T_n$ holds
  $$ \frac{p(\t)}{3n} - \frac 53 - n < d(\t) \le 512  (p(\t))^2 + 869
  p(\t) + 376\text{.\qed}$$ 
\end{theorem}

\section{Proof of Theorem~\ref{thm:crossing}}
\label{sec:proof}

\begin{proof}[Proof of the first part of Theorem~\ref{thm:crossing}]
  Let $\t\in\mathbb T_n$ be a polytopal triangulation of $S^3$ with $n$ 
  tetrahedra and let $L\subset  \t^1$ be a link formed by $k$ edges.
  Since $\t$ has at most $2n$ edges, it suffices to prove $\CR(L)< k^2$.  
  Since $\t$ is polytopal, it has a so-called 
  \emph{Schlegel diagram}~\cite{ziegler}, i.e., the simplicial
  complex obtained from $\t$ by removing a single tetrahedron 
  can be embedded into $\mathbb R^3$ so that any simplex is 
  euclidian. 
  In particular, the link $L\subset  \t^1$ is a union of $k$ straight line
  segments. It is possible to obtain a link diagram for $L$ by an orthogonal
  projection onto some plane. 
  Since an orthogonal  projection yields at most one crossing for any
  pair of edges, we obtain $\CR(L)< k^2$.
\end{proof}

In order to prove the second and third part of
Theorem~\ref{thm:crossing}, we formulate and prove the following
estimate for $\CR(L)$ in terms of $p(\t)$.
\begin{theorem}\label{thm:crosspoly}
  Let $\t$ be a triangulation of $S^3$, and let $L\subset \t^1$ be a
  link formed by $k$ edges of $\t$. Then 
  $$ \CR(L) < \big( k +  512 (p(\t))^2 + 869 p(\t) + 376\big)^2.$$
\end{theorem}
\begin{proof}
  By Theorem~\ref{thm:makepoly}, there is a polytopal triangulation
  $\tilde\t$ of $S^3$, obtained from $\t$ be a series of at most $512
  (p(\t))^2 + 869 p(\t) + 376$ expansions. 
  If a triangulation of $S^3$ is obtained from $\t$ by a single expansion, then
  its 1--skeleton contains a copy of $L$ formed by at most $k+1$
  edges. Thus, by induction on the number of expansions, there
  is a link $\tilde L\subset \tilde \t^1$ equivalent to $L$, formed by
  at most $k +  512 (p(\t))^2 + 869 p(\t) + 376$ edges. 
  Since $\tilde \t$ is polytopal, we obtain $\CR(L) = \CR(\tilde L) <
  \big( k +  512 (p(\t))^2 + 869 p(\t) + 376\big)^2$  as in the proof
  of the first part of Theorem~\ref{thm:crossing}.   
\end{proof}

\begin{proof}[Proof of the second and third part of
  Theorem~\ref{thm:crossing}] 
  Let  $\t\in\mathbb T_n$, and let $L\subset \t^1$ be a
  link. Since $L$ is formed by at most $2n$ edges,
  Theorems~\ref{thm:polybounds} and~\ref{thm:crosspoly} immediately
  yield 
  \begin{eqnarray*}
    \CR(L) &<& \big(25088 n^2  + 6085 n + 376\big)^2\\
    &<& 10^9 n^4
  \end{eqnarray*}
  provided $\t$ is shellable, and in general
  \begin{eqnarray*}
    \CR(L) &<& \Big(512\cdot 2^{400 n^2} + 869\cdot 2^{200 n^2} + 2n +
    376\Big)^2 \\
    &<& 2^{810 n^2}.
  \end{eqnarray*}    
\end{proof}


\begin{thebibliography}{}

\bibitem{armentrout}
  \newblock Armentrout, S.:
  \newblock Knots and shellable cell partitionings of $S^3$.
  \newblock {\em Illinois J. Math.} \textbf{38,} 347--365, (1994).

\bibitem{bruggesser}
  \newblock Bruggesser, H. and Mani, P.:
  \newblock Shellable decompositions of cells and spheres.
  \newblock {\em Math. Scand} \textbf{29,} 197--205 (1971).

\bibitem{ehrenborghachimori}
  \newblock Ehrenborg, R. and Hachimori, M.:
  \newblock Non-constructible complexes and the bridge index.
  \newblock Preprint 2000.

\bibitem{king1}
  \newblock King, S.:
  \newblock The size of triangulations supporting a given link.
  \newblock {\em Geometry \& Topology} \textbf{5,} 369--398 (2001).
  \newblock See~\cite{king2} for the slightly better bound that we
  use here.

\bibitem{king2}
  \newblock King, S.:
  \newblock {\em Polytopality of  triangulations}.
  \newblock PhD thesis,\\
  \texttt{http://www-irma.u-strasbg.fr/irma/publications/2001/01017.shtml.} 

\bibitem{lickorish}
  \newblock Lickorish, W. B. R.:
  \newblock Unshellable triangulations of spheres.
  \newblock {\em European J. Combinatorics} \textbf{12,} 527--530 (1991).

\bibitem{moise}
  \newblock Moise, E. E.:
  \newblock Affine structures in 3--manifolds.
  \newblock {\em Ann. Math.} \textbf{56,} 96--114 (1952).

\bibitem{schubert}
  \newblock Schubert, H.:
  \newblock \"Uber eine numerische Knoteninvariante.
  \newblock {\em Math. Z.} \textbf{61,} 245--288 (1954).

\bibitem{ziegler}
  \newblock Ziegler, G.:
  \newblock {\em Lectures on polytopes},
  \newblock Graduate Texts in Mathematics 152, 
  \newblock (Springer Verlag 1995).

\end{thebibliography}
\end{document}